\documentclass[11.5pt]{article}
\usepackage{mathrsfs}
\usepackage{color}
\usepackage{indentfirst}
\usepackage{amsmath,url,cases,supertabular}
\usepackage{amssymb,mathtools,multirow}
\usepackage{enumerate,makecell}
\usepackage{amsfonts}
\usepackage{amsthm}
\usepackage{graphicx,epstopdf}
\usepackage{color,verbatim}
\usepackage{mathrsfs,subeqnarray,subfigure}
\usepackage{listings,array}
\usepackage{booktabs}
\usepackage{dashrule}
\usepackage{arydshln}
\usepackage[vlined,ruled]{algorithm2e}
\usepackage{amsmath}
\usepackage{stmaryrd}
\usepackage{txfonts}
\usepackage{amsfonts}
\usepackage{amsfonts}
\usepackage{amssymb}
 \usepackage[table,xcdraw]{xcolor}
\RequirePackage{color,graphicx}
\def\rev#1{{\color{black}{{#1}}}}
\usepackage{cite}

\newcommand{\ba}{\noindent $\begin{array}}
\newcommand{\ea}{\end{array}$}
\newcommand{\be}{\begin{equation}}
\newcommand{\ee}{\end{equation}}
\newcommand{\bd}{\begin{displaymath}}
\newcommand{\ed}{\end{displaymath}}
\newcommand{\beq}{\begin{eqnarray*}}
\newcommand{\eeq}{\end{eqnarray*}}
\newcommand{\beqn}{\begin{eqnarray}}
\newcommand{\eeqn}{\end{eqnarray}}

\makeatletter

\newcommand{\Rmnum}[1]{\expandafter\@slowromancap\romannumeral #1@}
\makeatother

\usepackage{latexsym,lscape}

\newtheorem{theorem}{Theorem}[section]
\newtheorem{proposition}{Proposition}[section]
\newtheorem{lemma}{Lemma}[section]
\newtheorem{corollary}{Corollary}[section]

\newtheorem{remark}{Remark}[section]

\newtheorem{example}{Example}[section]

\newfont{\Bb}{msbm10 scaled\magstep1}

\def\sqr#1#2{{\vcenter{\hrule height .#2pt
      \hbox{\vrule width .#2pt height#1pt \kern#1pt\vrule width.#2pt}
                       \hrule height.#2pt}}}

\title{The Augmented Lagrangian Method Can Approximately Solve Convex Optimization  with Least Constraint Violation
}

\author{Yu-Hong Dai\footnote{LSEC, ICMSEC, AMSS, Chinese Academy of Sciences, Beijing 100190, China.  {\sl Email}: dyh@lsec.cc.ac.cn.
This author was supported by the Natural Science Foundation of China (Nos. 11991020, 12021001, 11631013, 11971372 and 11991021) and
the Strategic Priority Research Program of Chinese Academy of Sciences (No. XDA27000000).}
 \footnote{School of Mathematical Sciences, University of Chinese Academy of Sciences, Beijing 100049, China.}\quad and \quad Liwei Zhang \footnote{School of Mathematical Sciences, Dalian University of Technology, Dalian 116024, China.\ {\sl Email}: lwzhang@dlut.edu.cn. This author was supported by the Natural Science Foundation of China (Nos. 11971089 and 11731013) and partially supported by Dalian High-level Talent Innovation Project (No. 2020RD09).}}

\begin{document}

\maketitle
\begin{abstract}
There are many
important practical optimization problems whose feasible regions are not
known to be nonempty or not, and optimizers of the objective function with the least
constraint violation prefer to be found. A natural way for dealing with these problems
is to extend the nonlinear optimization problem as the one optimizing the objective
function over the set of points with the least constraint violation. This leads to the study of the  shifted problem. This paper focuses on the constrained convex optimization problem.
The sufficient condition for the closedness of the set of feasible shifts is presented and the continuity properties of the optimal value function and the solution mapping for the shifted problem are studied. Properties of the conjugate dual of the shifted problem are discussed through the relations between the dual function and the optimal value function. The solvability of the dual of the optimization problem with the least constraint violation is investigated. It is shown that, if the least violated shift is in the domain of the subdifferential of the optimal value function, then this dual problem has an unbounded solution set. Under this condition, the optimality conditions for the problem with the least constraint violation are established in term of the augmented Lagrangian. It is shown that the augmented Lagrangian method has the properties  that the sequence of shifts converges to the least violated shift and  the sequence of multipliers is  unbounded. Moreover, it is proved that the augmented Lagrangian method is able to find an approximate solution to the problem with the least constraint violation.

\vskip 12 true pt \noindent \textbf{Key words}: convex optimization, least constraint violation, augmented Lagrangian method, shifted problem, optimal value mapping, solution mapping, dual function, conjugate dual.
\end{abstract}


\section{Introduction}
\setcounter{equation}{0}
For a practical optimization problem, it is often unknown whether the feasible set of this problem is nonempty or not. In many cases, to design efficient numerical algorithms, infeasibility detection becomes an important issue when the constraint set is
 empty or the constraints are inconsistent. For nonlinear optimization, many numerical algorithms have been proposed to find infeasible stationary points; namely, stationary points for minimizing an infeasibility measure of the violated constraints.  Byrd, Curtis and Nocedal \cite{Byrd2010} presented a set of conditions to guarantee the superlinear convergence of their SQP algorithm to an infeasible stationary
point.  Burke, Curtis and Wang \cite{Burke2014} considered the general program
with equality and inequality constraints, and proved that their SQP method has
strong global convergence and rapid convergence to the KKT point, and has
superlinear/quadratic convergence to an infeasible stationary point. Recently, Dai, Liu and Sun \cite{DLSun2020}
proposed a  primal-dual interior-point method, which
can be superlinearly or quadratically convergent to the KKT point if the original problem is feasible, and can be superlinearly or quadratically convergent to the infeasible stationary point when the problem
is infeasible.

 These algorithms can find a stationary point of the infeasibility measure, which have nothing to do with the objective function of the problem. In practice, there are many problems that we need to find  minimizers of the objective function over the set of points with the least constraint violation. A natural way to deal with such problems is to  extend the constrained optimization problem as the one that optimizes the objective function over the set of points with least constraint violation. When the feasible region is nonempty, the set of points with least constraint violation coincides with the feasible region of the constrained optimization problem and hence the extended constrained optimization problem coincides with the original problem.

For the convex nonlinear optimization problem with the least constraint violation, Dai and Zhang \cite{DZhang2021} reformulated the problem as an optimization problem with a Lipschitz continuous equality (\cite{Clarke83}) or an MPCC problem (\cite{LPRalph96}) and derived necessary optimality conditions in terms of $L$-stationary point and $M$-stationary point. They also made efforts to construct a penalty function method for the problem with the least constraint violation and a smoothing algorithm for solving the resulted MPCC problem.
However, the penalty method is not an exact one, which will lead to computational difficulty when the penalty parameter is quite large and the smoothing method can only guarantee to converge to an $L$-stationary point  of the  problem with the least constraint violation. Thus it is worthwhile to reconsider both theoretical issues and algorithmic issues for the optimization problem with the least constraint violation. For example, the solvability of the problem, the characterization of optimal solutions, the characterization of the dual problem and what kind of algorithms can find an optimal solution or an approximate
solution. For the algorithmic issue, the work by  Chiche and  Gilbert \cite{ChCh2016}  gives us a valuable clue. They proved that the augmented Lagrangian method can deal with an infeasible convex quadratic optimization problem.  This stimulates us to consider the augmented Lagrangian method for dealing with a general convex optimization problem with the least constraint violation. Chiche and  Gilbert \cite{ChCh2016} also presented many theoretical results for the least shifted problem of the convex optimization problem, including the optimality conditions and differential characterization of the dual function and the unboundedness of the dual problem. All these results are quite helpful to us in studying the general convex optimization problem.

\rev{It should be pointed out that there are other interesting backgrounds similar to the convex optimization problem with the least constraint violation, one of which is the data-compatibility approach to constrained optimization proposed by Censor,  Zaknoon and  Zaslavski \cite{Censor2021}.}

The augmented Lagrangian  method is a famous method for solving constrained optimization problems. It was proposed  by Hestenes
\cite{H69} and Powell \cite{P72} for solving  optimization problems with equality constraints  and was generalized by
Rockafellar  \cite{R73a} to
  optimization problems with both equality and inequality constraints.
    Rockafellar \cite{R73a}  demonstrated  a
saddle point theorem for  convex  optimization in terms of the augmented Lagrangian,  and Rockafellar  \cite{R73b}
established
the global convergence of the augmented Lagrangian  method for convex optimization with any positive penalty parameter.
It was observed by Rockafellar \cite{Rockafellar76b} that, for a convex optimization problem with inequality constraints, the augmented Lagrangian  method  is essentially the proximal point method (for maximal monotone operator inclusion proposed in \cite{Rockafellar76a}) applied to the maximal monotone operator inclusion expressing the optimality conditions of the dual problem.  Many publications  studied the rate of convergence of  the augmented Lagrangian method for solving various feasible non-convex optimization problems, see references  \cite{P72},  \cite[Chapter 3]{B82}, \cite{CGToint91}, \cite{Contesse-Becker93}, \cite{IKunisch90} and \cite{SSZhang2008}.

In this paper, we shall study the augmented Lagrangian method for dealing with the convex optimization problem with the least constraint violation when the convex optimization problem is infeasible.

Consider the following constrained optimization problem,
\begin{equation}\label{eq:p1}
({\rm P})\quad \quad
\begin{array}{cl}
\min & f(x)\\
{\rm s.t.}& g(x)\in K,
\end{array}
\end{equation}
where $f: \Re^n \rightarrow \Re$, $g:\Re^n \rightarrow {\cal Y}$, $K \subset {\cal Y}$ is a nonempty  closed convex set and ${\cal Y}$ is a finite dimensional Hilbert space.
 The optimization model (\ref{eq:p1}) covers a lot of optimization problems:
\begin{itemize}
\item $K=\{0_q\}\times \Re^p_+$:\quad  nonlinear optimization problem;
\item $K=\{0_q\}\times \Re^p_+\times \mathcal{S}^m_+$:\quad  semidefinite optimization problem;
\item $K=\{0_q\}\times \Re^p_+\times [Q_{m_1+1}\times \cdots\times Q_{m_J+1}]$:\quad  second-order conic optimization problem;
\item $K=\{0_q\}\times \Re^p_+\times {\rm epi}\, \|\cdot\|_2$:\quad spectral norm induced conic optimization problem;
\item $K=\{0_q\}\times \Re^p_+\times {\rm epi}\, \|\cdot\|_*$:\quad nuclear norm induced conic optimization problem;
\item $K=\{0_q\}\times \Re^p_+\times {\rm epi}\, \|\cdot\|_k$:\quad $k$-norm induced conic optimization problem.
\end{itemize}
Therefore, we can deal with at least the above optimization problems with the least constraint violation if we are able to establish a theoretical framework and an algorithmic analysis for the general model (\ref{eq:p1}). For simplicity, we only consider the convex optimization problem.  We say that Problem (\ref{eq:p1}) is convex if $f$ is a convex function and the set valued mapping
\begin{equation}\label{eq:n0}
G_K(x):=-g(x)+K
\end{equation}
is graph-convex; namely, its graph
$$
{\rm gph}\,G_K=\left\{(x,y)\in\Re^n \times {\cal Y}:y \in -g(x)+K\right\}
$$
is a convex set in $\Re^n \times {\cal Y}$. In this paper, we assume that Problem  (\ref{eq:p1}) is convex and $f$ and $g$ are continuously differentiable.

 By introducing an artificial vector $y\in {\cal Y}$, Problem (\ref{eq:p1}) can equivalently be expressed as
 \begin{equation}\label{eq:p2}
\begin{array}{cl}
\min & f(x)\\
{\rm s.t.}& g(x)=y,\\[0pt]
& y\in K.
\end{array}
\end{equation}
In Section \ref{sec:ALM}, we shall present the augmented Lagrangian method based on the augmented Lagrangian of Problem (\ref{eq:p2}), which will bring us some convenience for analysis.

The rest of this paper is organized as follows. In Section 2, we consider the shifted problem  of Problem (\ref{eq:p1}) and discuss properties of the optimal value and the set of optimal solutions. We provide conditions for the closedness of the set of feasible lifts, the conditions for the existence of solution as well as the continuity properties of the optimal value and the set of optimal solutions depending on the shifts. In Section 3, the conjugate dual of the shifted problem is studied. From the classical duality theory, we obtain a characterization of the solution set of the dual in terms of the dual function and the optimal value function. Especially, the equivalence of the subdifferentiability of the dual function and the existence of an augmented Lagrangian problem are established. This is crucial for the analysis of the augmented Lagrangian method. Section 4 is the central part of this paper. The optimality conditions (both necessary and sufficient) for the convex optimization problem with the least constraint violation are established in terms of the augmented Lagrangian, which suggest us that the augmented Lagrangian method is a suitable method for solving this problem. The classical augmented Lagrangian method is analyzed and its shown that the dual of the least constraint violated problem has an unbounded solution set if it has a solution or if the least violated shift is in the domain of the subdifferential of the optimal value function. Under this condition, it is  shown that the augmented Lagrangian method has the properties  that the sequence of shifts converges to the least violated shift and it is able to find an approximate solution. Section 5 gives some discussions about possible topics on optimization with the least constraint violation.

\section{The shifted problem}
\setcounter{equation}{0}
In this section, we shall define the shifted problem associated with Problem  (\ref{eq:p1}) and the set of feasible shifts, and discuss when the optimal value function is lower semi-continuous and when the set of feasible shifts is closed.

For a given $s \in {\cal Y}$, the shifted problem is defined as
 \begin{equation}\label{eq:p3}
 {\rm P}(s)\quad \quad
 \begin{array}{ll}
  \min & f(x)\\
{\rm s.t.}& g(x)+s\in K.
\end{array}
\end{equation}
By  introducing an artificial vector $y\in {\cal Y}$, Problem (\ref{eq:p3}) can equivalently be expressed as
 \begin{equation}\label{eq:p3y}
  \begin{array}{ll}
\min & f(x)\\
{\rm s.t.}& g(x)+s=y,\\[0pt]
& y\in K.
\end{array}
\end{equation}
Here we call $s$ as a shift. The set of feasible shifts, denoted as ${\cal S}$, is defined by
\begin{equation}\label{eq:n1}
{\cal S}:=\{s \in {\cal Y}: \mbox{there exists some } x \in \Re^n \mbox{ such that } g(x)+s \in K\}.
\end{equation}
Let
$$\Phi(s)=\{x\in \Re^n: g(x)+s\in K\}
$$
be the feasible set of Problem ${\rm P}(s)$. Denote by $\nu (s)$ and $X(s)$ the optimal value and the set of optimal solutions of Problem ${\rm P}(s)$, respectively; namely,
\begin{equation}\label{eq:val}
\nu(s)=\inf_x\,\{ f(x): x \in \Phi(s)\} \quad \mbox{and}\quad X(s)={\rm argmin}\,\{ f(x): x \in \Phi(s)\}.
\end{equation}
 Obviously, if $s \not \in {\cal S}$, then $\Phi(s)=\emptyset$, $X(s)=\emptyset$  and $\nu (s)=+\infty$. If $f$ is a proper lower semi-continuous function over $\Re^n$, then
 \begin{equation}\label{eq:S2}
 {\rm dom }\, \nu ={\cal S}.
 \end{equation}
Noting that
$$
{\rm gph}\,G_K=\{(x,s) \in \Re^n\times {\cal Y}: s \in G_K(x)\}=\{(x,s) \in \Re^n\times {\cal Y}:g(x)+s \in K\},
$$
we have
\begin{equation}\label{eq:n2}
{\cal S}=\Pi_{{\cal Y}}({\rm gph}\,G_K),
\end{equation}
where $\Pi_{{\cal Y}}:\Re^n\times {\cal Y}\rightarrow {\cal Y}$ is the following projection mapping
$$
\Pi_{{\cal Y}}(x,s)=s\rev{,} \quad \forall (x,s) \in \Re^n\times {\cal Y}.
$$
The closedness of ${\cal S}$ is crucial for the existence of the least violated shift since the shift with the smallest norm is just the projection of zero onto ${\cal S}$. The following result gives a sufficient condition ensuring that ${\cal S}$ is closed.
\begin{lemma}\label{lem-1}
Let $f$ be a continuous convex function and $g$ be a continuous mapping such that $G_K$ is a graph-convex set-valued mapping.
If
\begin{equation}\label{cond-1}
0 \in G_K^{\infty}(h_x) \Longrightarrow h_x=0,
\end{equation}
then ${\cal S}$  is a closed set in ${\cal Y}$.
\end{lemma}
{\bf Proof}. In view of the equality ${\cal S}=\Pi_{{\cal Y}}({\rm gph}\,G_K)$, we have from Theorem 3.10 of \cite{RW98} that a sufficient condition for the closedness of ${\cal S}$ is
\begin{equation}\label{eq:h2h}
\Pi_{\cal Y}^{-1}(0) \cap ({\rm gph}\,G_K)^{\infty}=\{0\}\subset \Re^n \times {\cal Y}.
\end{equation}
From the definition of ${\rm gph}\,G_K^{\infty}$, we know that (\ref{eq:h2h}) is equivalent to
\begin{equation}\label{eq:h2a}
\Pi_{\cal Y}^{-1}(0) \cap {\rm gph}\,G_K^{\infty}=\{0\}\subset \Re^n \times {\cal Y}.
\end{equation}
Since $\Pi_{\cal Y}^{-1}(0)=\Re^n \times \{0\}$, (\ref{eq:h2a}) is equivalent to
$$
\rev{\left(\Re^n \times \{0\} \right)}\cap {\rm gph}\,G_K^{\infty}=\{0\}\subset \Re^n \times {\cal Y},
$$
 which is just (\ref{cond-1}).  The proof is completed. \hfill $\Box$
\begin{proposition}\label{prop-1}
Let $f$ be a continuous convex function and $g$ be a continuous mapping such that $G_K$ is a graph-convex set-valued mapping.
If
\begin{equation}\label{cond-1a}
g^{\infty} (h_x) \in K^{\infty} \Longrightarrow h_x=0,
\end{equation}
then ${\cal S}$ is closed, where $K^{\infty}$ is the recession cone of $K$ in the sense of convex analysis and $g^{\infty}$ is the horizon mapping defined by
$$
g^{\infty}(h):=\{u\in {\cal Y}: \exists\, t_k \rightarrow 0, \exists\, v^k \in \Re^n \mbox{ such that } t_k (v^k, g(v^k))\rightarrow (h,u)\}.
$$
\end{proposition}
{\bf Proof}.  Since $G_K$ is graph-convex, ${\rm ghp}\, G_K$ is a convex set in $\Re^n \times {\cal Y}$.
From the definition of $g^{\infty}$, for some $(x^0,s^0)$ satisfying $g(x^0)+s^0 \in K$, we have that
\begin{equation}\label{eq:h1}
\begin{array}{rcl}
({\rm ghp}\, G_K)^{\infty}& =&\left\{(h, u)\in \Re^n \times {\cal Y}: (x^0,s^0)+t (h,u) \in {\rm ghp}\, G_K,\,\forall t \geq 0\right\}\\[4pt]
&=&\left\{(h, u)\in \Re^n \times {\cal Y}: g(x^0+t h)+(s^0+t u) \in K,\forall t \geq 0\right\}\\[4pt]
&=&\left\{(h, u)\in \Re^n \times {\cal Y}: g(x^0)+s^0+t\displaystyle \left[u+\displaystyle
\frac{g(x^0+t h)-g(x^0)}{t}\right] \in K,\forall t \geq 0\right\}\\[4pt]
& \subseteq & \{(h, u)\in \Re^n \times {\cal Y}:u +g^{\infty}(h)\in K^{\infty}\}.
\end{array}
\end{equation}
Since the condition (\ref{cond-1a}) gives
\begin{equation}\label{eq:h2}
\Pi_{\cal Y}^{-1}(0) \cap \{(h, u)\in \Re^n \times {\cal Y}:u +g^{\infty}(h)\in K^{\infty}\} =\{0\}\subset \Re^n \times {\cal Y},
\end{equation}
we obtain from (\ref{eq:h1}) that
$$
\Pi_{\cal Y}^{-1}(0) \cap ({\rm gph}\,G_K)^{\infty}=\{0\}\subset \Re^n \times {\cal Y},
$$
which indicates the truth of (\ref{cond-1}).  The result follows from Lemma \ref{lem-1}. The proof is completed. \hfill $\Box$
\begin{corollary}\label{cor-1}
Let $K=\{0_q\}\times \Re^p_-$ and $g(x)=(a_1^Tx-b_1,\ldots, a_q^Tx-b_q, g_{q+1}(x),\ldots,g_{q+p}(x))^T$ with each $g_i$ being a lower semicontinuous convex function for $i=q+1,\ldots,q+p$.
If
$$
\left.
\begin{array}{l}
a_j^Tw=0, \ j =1,\ldots, q\\
g_i^{\infty}(w)\leq 0,\ i=1,\ldots, p
\end{array}
\right\}
\Longrightarrow w=0,
$$
then ${\cal S}$ is closed.
\end{corollary}
Define $\phi: \Re^n \times {\cal Y}  \rightarrow \overline{\Re}$ by
\begin{equation}\label{eq:phi}
\phi(x,s)=f(x)+\delta_K(g(x)+s),
\end{equation}
where $\delta_K$ is the indicator function
$$
\delta_K (y)=\left
\{
\begin{array}{ll}
0, & y \in K;\\[3pt]
+\infty,\quad & y \notin K.
\end{array}
\right.
$$
\begin{remark}\label{remark-1}
It follows from Proposition 2.162 of \cite{BS00}
that $(x,s) \rightarrow \delta_K(g(x)+s)$ is a convex function when $G_K$ is a graph-convex set-valued mapping.
\end{remark}
Now, we discuss, for $ s\in {\cal S}$,  when Problem ${\rm P}(s)$ has a nonempty set of  solutions and the continuity properties of $X(s)$ and $\nu(s)$. This requires the uniform level boundedness condition from Definition 1.16 of \cite{RW98}.
\begin{proposition}\label{prop:ulb}
Let $f$ be a continuous convex function and $g$ be a continuous mapping such that $G_K$ is a graph-convex set-valued mapping.
Then $\phi (x,s)$ is level-bounded in $x$ locally uniformly in $s$ in the sense of Definition 1.16 of \cite{RW98} if and only if
\begin{equation}\label{eq:inpc}
\left.
\begin{array}{l}
f^{\infty}(h_x)\leq 0\\[3pt]
0 \in G_K^{\infty}(h_x)
\end{array}
\right
\} \Longrightarrow h_x=0.
 \end{equation}
\end{proposition}
{\bf Proof}. From the definition of $\phi$, it is easy to check
\begin{equation}\label{eq:h3}
\phi(x,s)=f(x)+\delta_{\mbox{ghp }G_K}(x,s).
\end{equation}
In view of Remark \ref{remark-1}, we know that $\phi(x,s)$ is a proper lower semi-continuous convex function.
It follows from Theorem 3.31 of \cite{RW98} that $\phi(x,s)$ is level-bounded in $x$ locally uniformly in $s$ if and only if
\begin{equation}\label{eq:hc}
\phi^{\infty}(h_x,0)>0,\ \ \forall h_x\ne 0\in \Re^n.
\end{equation}
Since $f$ is continuous convex, ${\rm dom}\, f =\Re^n$, we have
$$
{\rm dom }\, f \times {\cal Y} \cap {\rm gph}\,G_K \ne \emptyset.
$$
Then, from Exercise 3.29 of \cite{RW98} and the identity that $\delta_{C}^{\infty}=\delta_{C^{\infty}}$, we obtain
$$
\phi^{\infty}(h_x,h_s)=f^{\infty}(h_x)+\delta_{({\rm ghp}\,G_K)^{\infty}}(h_x,h_s)=f^{\infty}(h_x)+\delta_{{\rm ghp}\,G_K^{\infty}}(h_x,h_s).
$$
Thus, the condition (\ref{eq:hc}) is equivalent to
$$
f^{\infty}(h_x)+\delta_{{\rm ghp}\,G_K^{\infty}}(h_x,0)>0, \ \ \forall  h_x\ne 0\in \Re^n,
$$
or
\begin{equation}\label{eq:hcb}
f^{\infty}(h_x)+\delta_{{\rm ghp}\,G_K^{\infty}}(h_x,0)\leq 0 \Longrightarrow  h_x= 0\in \Re^n.
\end{equation}
Obviously, the condition (\ref{eq:hcb}) is the same as the condition (\ref{eq:inpc}). The proof is completed. \hfill $\Box$
\begin{proposition}\label{prop:nulsc}
Let $f$ be a continuous convex function and $g$ be a continuous mapping such that $G_K$ is a graph-convex set-valued mapping.
Assume that the condition (\ref{eq:inpc}) is satisfied.
 Then
\begin{itemize}
\item[(a)] The
function $\nu$ is proper and lower semi-continuous convex on ${\cal Y}$, and for each $s \in
\mbox{dom}\, \nu$ the set $X(s)$ is nonempty, compact and convex, whereas
$X(s) =\emptyset$ when $s \notin \mbox{dom}\, \nu$.
\item[(b)] The mapping $X$, which is compact-valued with $\mbox{dom
}X=\mbox{dom }\nu$, is outer semi-continuous with respect to $\nu$-attentive convergence
$\stackrel{\nu}\rightarrow$.
  \item[(c)] The set-valued mapping $X$ is locally bounded and outer semi-continuous relative to the set $\mbox{int}(\mbox{dom
}X)= \mbox{int}(\mbox{dom }\nu)$.
\end{itemize}
\end{proposition}
{\bf Proof}. It follows from Proposition \ref{prop:ulb} that the condition (\ref{eq:inpc}) is equivalent to the condition that
$\phi (x,s)$ is level-bounded in $x$ locally uniformly in $s$ in the sense of Definition 1.16 of \cite{RW98}. The results follow from Theorems
1.17 and 7.41 of \cite{RW98}. \hfill $\Box$
\begin{remark}\label{remark-2}
Property (b) of the above proposition means that, if
$s^{k} \in X(s^{k})$ and if $s^{k} \rightarrow \hat s \in
\mbox{dom}\, \nu$ in such a way that $\nu(s^{k}) \rightarrow
\nu(\hat s)$ (as when $\nu$ is continuous at $\hat s$ relative to a
set $U$ containing $\hat{s}$ and $s^{k}$), then the sequence
$\{x^{k}\}_{k \in \textbf{N}}$ is bounded, and all its cluster
points lie in $X(\hat s)$.
\end{remark}
\begin{remark}\label{remark3}
Noting that ${\cal S}={\rm dom}\, \nu$, we have the following observations.
\begin{itemize}
\item[{\rm (i)}]
The first part of (c)  in  the above proposition means that, for any $s \in {\cal S}$, there exists $\varepsilon>0$ with a nonempty compact set $B \subset {\cal Y}$ such that
$$
X(s') \subset B\rev{, \quad \forall} s' \in {\rm dom}\, \nu \cap \textbf{B}(s,\varepsilon).
$$
\item[{\rm (ii)}]If $\mbox{int }{\cal S} \ne \emptyset$, then for any $s\in \mbox{int }{\cal S}$, $\nu$ is continuous at $s$
and
$$
\displaystyle \limsup_{s' \rightarrow s} X(s')=X(s).
$$
\item[{\rm (iii)}]If $\mbox{int }{\cal S} \ne \emptyset$, then, from the second part of (c) in  the above proposition, for any $s\in {\cal S}$, $\nu$ is lower semi-continuous at $s$ and
$$
\displaystyle \limsup_{s' \stackrel{{\rm int}\,{\cal S}}\longrightarrow s} X(s')=X(s).
$$
Since, for $s(t)=(1-t)s+t s_0$, where $s_0 \in {\rm int}\,{\cal S}$, one has
$$
s(t) \in  {\rm int}\,{\cal S}, \quad   \quad \mbox{    for   } \ t \in (0,1) \quad \mbox{    and     } \quad \lim_{t\rightarrow 0_+} s(t)=s,
$$
so that
$$
\displaystyle \limsup_{t \searrow 0} X(s(t))\subseteq X(s).
$$
\item[{\rm (iii)}]If $\mbox{int }{\cal S} =\emptyset$, in this case ${\rm ri}\, {\cal S}\ne \emptyset$, then one only has, from  (a)  in  the above proposition, that
$$
\displaystyle \limsup_{s' \stackrel{\nu} \rightarrow s} X(s')=X(s).
$$
\item [{\rm (iv)}] If $\mbox{int }{\cal S} =\emptyset$, in this case ${\rm ri}\, {\cal S}\ne \emptyset$, then for any $s\in \mbox{ri }{\cal S}$, it follows from Theorem 10.1 of \cite{Rock70}\footnote{
Theorem 10.1. A convex function $f$ on $\Re^n$ is continuous relative to any
relatively open convex set $C$ in its effective domain, in particular relative to
${\rm ri}\, ({\rm dom}\, f)$.} that $\nu$ is continuous at $s$ relative to
${\rm ri}\, ({\rm dom}\, \nu)$ so that
$$
\displaystyle \limsup_{s' \stackrel{{\rm ri}\,{\cal S}}\longrightarrow s} X(s')=X(s).
$$
For $s(t)=(1-t)s+t s_0$, where $s_0 \in {\rm ri}\,{\cal S}$, one has
$$
s(t) \in  {\rm ri}\,{\cal S}  \quad \mbox{    for   } \quad t \in (0,1) \quad \mbox{    and     } \quad \lim_{t\rightarrow 0_+} s(t)=s,
$$
so that
$$
\displaystyle \limsup_{t \searrow 0} X(s(t))\subseteq X(s).
$$
\end{itemize}
\end{remark}

Now we pay a special attention to convex quadratic optimization.

\begin{example}\label{ex:qp}
Consider the following convex quadratic optimization problem
\begin{equation}\label{eq:pqp}
\begin{array}{cl}
\min & f(x):=c^Tx+\displaystyle \frac{1}{2}x^TGx\\[3pt]
{\rm s.t.} & Ex=d,\\[3pt]
& Ax \leq b,
\end{array}
\end{equation}
where $c\in \Re^n$, $G \in \Re^{n\times n}$ is a positive semidefinite symmetric matrix, $E \in \Re^{q\times n}$, $d \in \Re^q$,
$A \in \Re^{p\times n}$ and $b \in \Re^p$. Denote $s=(s_e,s_a)\in \Re^{q+p}$ with $s_e \in \Re^q$ and $s_a \in \Re^p$.

 For a given $s \in \Re^{q+p}$, the shifted problem becomes
 \begin{equation}\label{eq:p3qp}
 {\rm QP}(s)\quad \quad
 \begin{array}{ll}
\min & f(x):=c^Tx+\displaystyle \frac{1}{2}x^TGx\\[3pt]
{\rm s.t.} & Ex+s_e=d,\\[3pt]
& Ax +s_a\leq b.
\end{array}
\end{equation}
The set of feasible shifts is given by
\begin{equation}\label{eq:n1-b}
{\cal S}_{QP}:=\{s \in {\cal Y}: \mbox{there exists some } x \in \Re^n \mbox{ such that } Ex+s_e=d,\, Ax +s_a\leq b\}.
\end{equation}
Let
$\nu (s)$ and $X(s)$ be the optimal value and the set of optimal solutions of Problem ${\rm QP}(s)$, respectively; namely,
$$
\begin{array}{rcl}
\nu_{QP}(s)&=&\inf_x\{c^Tx+\displaystyle \frac{1}{2}x^TGx: Ex+s_e=d,\, Ax +s_a\leq b\}, \\
 X_{QP}(s)&=&{\rm argmin}\,\{c^Tx+\displaystyle \frac{1}{2}x^TGx: Ex+s_e=d,\, Ax +s_a\leq b\}.
\end{array}
$$
 The set-valued mapping $G_K: \Re^n \rightrightarrows \Re^{q+p}$ is given by
$$
G_K(x)=\left[
\begin{array}{l}
-Ex+d\\[4pt]
-Ax+b
\end{array}
\right]+\{0_q\} \times \Re^p_-.
$$
Then we obtain
$$
{\rm ghp}\, G_K=\{(x,s) \in \Re^n \times \Re^{q+p}:Ex+s_e=d,\, Ax+s_a\leq b\}
$$
and
$$
\left({\rm ghp}\, G_K\right)^{\infty}=\{(h_x,h_s) \in \Re^n \times \Re^{q+p}:Eh_x+h_{s_e}=0,\,Ah_x+h_{s_a}\leq 0\}\rev{,}
$$
so that
\begin{equation}\label{eq:qphr}
G_K^{\infty}(h_x)=\{h_s\in \Re^{q+p}:Eh_x+h_{s_e}=0,\,Ah_x+h_{s_a}\leq 0\}.
\end{equation}
\begin{itemize}
\item[{\rm (a)}]
In view of Lemma \ref{lem-1}, we have that, if
\begin{equation}\label{eq:condpq}
\left.
\begin{array}{l}
Eh_x=0\\[3pt]
Ah_x \leq 0
\end{array}
\right
\} \Longrightarrow h_x=0,
  \end{equation}
then ${\cal S}_{QP}$ is closed;
\item[{\rm (b)}]
It follows from page 89 of \cite{RW98} that
$$
f^{\infty}(h_x)=c^Th_x+\delta (h_x\,|\, Gh_x=0).
$$
Thus the condition (\ref{eq:inpc}) in this case becomes
\begin{equation}\label{eq:qpc}
\left.
\begin{array}{l}
c^Th_x\leq 0,\
Gh_x=0\\[3pt]
Eh_x=0,\
Ah_x \leq 0
\end{array}
\right
\} \Longrightarrow h_x=0.
  \end{equation}
  It follows from Proposition \ref{prop:nulsc} that, if the condition (\ref{eq:qpc}) is satisfied, then
\begin{itemize}
\item[{\rm (b1)}] The
function $\nu_{QP}$ is proper and lower semi-continuous convex on $\Re^{q+p}$, and for each $s \in
{\cal S}_{QP}$ the set $X_{QP}(s)$ is nonempty, compact and convex.
\item[{\rm (b2)}] The mapping $X_{QP}$, which is compact-valued with $\mbox{dom
}X_{QP}={\cal S}_{QP}$, is outer semi-continuous with respect to $\nu_{QP}$-attentive convergence
$\stackrel{\nu_{QP}}\rightarrow$;
 \end{itemize}
\item[{\rm (c)}] Moreover, if there are no equality constraints in Problem (\ref{eq:pqp}), then $\mbox{int }{\cal S}_{QP}\ne \emptyset$. In this case, the condition (\ref{eq:qpc}) is reduced to the following condition
 \begin{equation}\label{eq:qpcne}
\left.
\begin{array}{l}
c^Th_x\leq 0\\[3pt]
Gh_x=0\\[3pt]
Ah_x \leq 0
\end{array}
\right
\} \Longrightarrow h_x=0.
  \end{equation}
Thus, if (\ref{eq:qpcne}) is satisfied, then the set-valued mapping $X_{QP}$ is locally bounded and outer semi-continuous  relative to the set $\mbox{int } {\cal S}_{QP}$;
\item[{\rm (d)}]One has from Lemma 2.2 of \cite{ChCh2016}, for $s \in {\cal S}_{QP}$, that Problem ${\rm P}_{QP}(s)$ is unbounded if and only if  there exists $h_x \in \Re^n$ such that
\begin{equation}\label{eq:qpubd}
\begin{array}{l}
c^Th_x<0,\ \
Gh_x=0,\\[3pt]
Eh_x=0,\ \
Ah_x \leq 0.
\end{array}
  \end{equation}
\end{itemize}
  \end{example}
  Obviously, if ${\cal Y}\ni 0\in {\cal S}$, then Problem (\ref{eq:p1}) is feasible. Otherwise, we define the smallest norm shift, denoted by $\bar s$,  as the projection of $0 \in {\cal Y}$ on ${\cal S}$:
  \begin{equation}\label{eq:snshiftorig}
    \bar s=\mbox{argmin} \left\{\displaystyle \frac{1}{2}\|s\|^2: s \in {\cal S}\right\}.
  \end{equation}
  If ${\cal S}$ is closed, then $\bar s$ can be achieved; namely, $\bar s \in {\cal S}$. In this case, the optimization problem with the least constraint violation is expressed as follows
   \begin{equation}\label{eq:1.12orig}
   {\rm P}(\bar s)\quad \quad
 \begin{array}{ll}
  \min & f(x)\\
{\rm s.t.}& g(x)+\bar s\in K.
\end{array}
\end{equation}

\section{The dual of the shifted problem}
\setcounter{equation}{0}
This section will discuss properties of the dual problem through the relations between the dual function and the optimal value function.
As shown in the next section, the analysis here will help us to understand the behaviour of the optimization problem with least constraint violation.

The Lagrangian of Problem (\ref{eq:p2}), denoted by $l:\Re^n \times {\cal Y}\times {\cal Y}$, is defined by
\begin{equation}\label{eq:2.5}
l(x,y,\lambda)=f(x)+\langle \lambda, g(x)-y\rangle.
\end{equation}
The augmented Lagrangian function of Problem (\ref{eq:p2}), denoted by $l_r:\Re^n \times {\cal Y}\times {\cal Y}$, is defined by
\begin{equation}\label{eq:2.5agm}
l_r(x,y,\lambda)=f(x)+\langle \lambda, g(x)-y\rangle+\displaystyle \frac{r}{2}\|g(x)-y\|^2.
\end{equation}
The dual function $\theta:{\cal Y} \rightarrow \overline \Re$ associated with Problem (\ref{eq:p2}) is
\begin{equation}\label{eq:2.6}
\theta(\lambda):=-\inf_{x\in \Re^n, y \in K} l(x,y,\lambda).
\end{equation}
It is easy to check that
\begin{equation}\label{eq:2.6a}
\theta(\lambda):=-\inf_{x\in \Re^n} \left\{f(x)+\langle\lambda, g(x)\rangle-\delta^*(\lambda \,|\, K)\right\},
\end{equation}
where $\delta^*(\lambda\,|\,K)$ is the supporting function of $K$ at $\lambda$:
$$
\delta^*(\lambda\,|\,K)=\sup_{y \in K} \langle y, \lambda \rangle.
$$
The function $\theta$ is a lower semi-continuous convex function, but does not take value $-\infty$. Therefore, we have
\begin{equation}\label{eq:2.7}
\theta \mbox{ is a proper lower semi-continuous convex function }\Longleftrightarrow {\rm dom}\,\theta \ne \emptyset.
\end{equation}

\begin{proposition}\label{prop:2.7}
Suppose ${\rm dom}\, \theta\ne \emptyset$ and ${\cal S}$ is closed. Then the following two properties are equivalent:
\begin{itemize}
\item[{\rm (i)}] Problem (\ref{eq:p1}) is feasible;
\item[{\rm (ii)}] The dual function $\theta$ is bounded below.
\end{itemize}
\end{proposition}
{\bf Proof}. (i)$\Longrightarrow$(ii). Assume that Problem (\ref{eq:p1}) be feasible. There is some $x_0$ such that $y_0=g(x_0)\in K$. It follows from the definition of $\theta$ that, for any $\lambda \in {\cal Y}$,
$$
\theta (\lambda)=-\inf_{x\in \Re^n, y \in K} l(x,y,\lambda)\geq -l(x_0,y_0,\lambda)=-f(x_0),
$$
implying that $\theta$ is bounded below by $-f(x_0) \in \Re$.

(ii)$\Longrightarrow$(i). Since ${\rm dom}\, \theta\ne \emptyset$, there exists $\lambda \in {\cal Y}$ such that $\theta (\lambda)\in \Re$. One the other hand, since
$$
 \bar s=\mbox{argmin} \left\{\psi(s):=\displaystyle \frac{1}{2}\|s\|^2: s \in {\cal S}\right\}
$$
and ${\cal S}$ is an nonempty closed convex set, we have that
$$
{\rm D}\psi (\bar s)(s-\bar s)\geq 0 \rev{,}\quad \forall s \in {\cal S},
$$
or
$$
\langle s, \bar s\rangle \geq \|\bar s\|^2 \rev{,}\quad \forall s \in {\cal S}.
$$
Then for all $(x,y) \in \Re^n \times K$,
$$
\langle y- g(x),\bar s \rangle\geq \|\bar s\|^2.
$$
Then, for all $t\geq 0$,
$$
\begin{array}{rcl}
\theta (\lambda-t\bar s) &=&-\displaystyle\inf_{x \in \Re^n, y\in K} [f(x)+\langle \lambda-t\bar s, g(x)-y \rangle]\\[8pt]
&=&-\displaystyle\inf_{x \in \Re^n, y\in K} [f(x)+\langle \lambda,  g(x)-y \rangle-t\, \langle \bar s, g(x)-y \rangle]\\[8pt]
&\leq&-\displaystyle\inf_{x \in \Re^n, y\in K} [f(x)+\langle \lambda,  g(x)-y \rangle-t\, \| \bar s\|^2]\\[8pt]
&=& \theta (\lambda)-t\, \|\bar s\|^2.
\end{array}
$$
Since $\theta$ is bounded below, we must have that $\bar s=0$; namely, Problem (\ref{eq:p1}) is feasible. \hfill $\Box$
\begin{proposition}\label{prop:2.8}
For $\nu$ defined by
(\ref{eq:val})
and $\theta$ defined by (\ref{eq:2.6}). The following properties are satisfied:
\begin{itemize}
\item[{\rm (a)}]$
\theta (\lambda)=\nu^*(\lambda)$ and $\nu^{**}(s)=\theta^*(s);
$
\item[{\rm (b)}]If $\nu$ is a proper function (namely $\nu(s)>-\infty, \forall s \in {\cal S}$ and there exists a vector $\hat s \in {\cal S}$ such that $\nu (\hat s)<+\infty$), then $\theta$ and $\nu^{**}$ are proper lower semi-continuous convex functions (indicating that ${\rm dom}\,\theta\ne \emptyset$ and ${\rm dom}\,\nu^{**}\ne \emptyset$);
 \item[{\rm (c)}] If $\nu$ is a proper lower semi-continuous function, then for any $s \in {\cal Y}$,
$$
\nu(s)=\nu^{**}(s)=\theta^*(s).
$$
\end{itemize}
\end{proposition}
{\bf Proof}. Define
$$
{\cal L}(x,\lambda)=f(x)+\langle \lambda, g(x)\rangle.
$$
Then
\begin{equation}\label{eq:aa}
\begin{array}{rcl}
\nu^*(\lambda)& =&\displaystyle\sup_{s\in {\cal Y}} \left \{ \langle \lambda, s\rangle -\nu (s) \right\}\\[4pt]
&=&\displaystyle\sup_{s\in {\cal Y}} \left \{ \langle \lambda, s\rangle -\displaystyle \inf_x\,[f(x)+ \delta_K(g(x)+s)] \right\}\\[4pt]
&=&\displaystyle\sup_{s\in {\cal Y}} \left \{ \langle \lambda, u-g(x)\rangle -\displaystyle \inf_x\,[f(x)+ \delta_K(u)] \right\}\\[4pt]
&=&\displaystyle\sup_{x\in \Re^n, u\in {\cal Y}}\left\{-f(x)-\langle \lambda, g(x)\rangle+\langle \lambda, u\rangle-\delta_K(u)\right\}\\[4pt]
&=&\delta^*_K(\lambda)-{\displaystyle \inf_{x\in \Re^n}}\, {\cal L}(x,\lambda).
\end{array}
\end{equation}
From the definition of $\theta (\lambda)$, we have that
\begin{equation}\label{eq:bb}
\begin{array}{rcl}
\theta (\lambda)&=&-\displaystyle \inf_{x\in \Re^n, y \in K} l(x,y,\lambda)\\[8pt]
&=&-\displaystyle\inf_{x\in \Re^n, y \in K} [f(x)+\langle \lambda, g(x)-y \rangle]\\[8pt]
&=&\displaystyle\sup_{x\in \Re^n, y \in K} [-f(x)-\langle \lambda, g(x)\rangle+\langle \lambda, y \rangle]\\[12pt]
&=&\displaystyle\sup_{x\in \Re^n} \left[-{\cal L}(x,\lambda)+\sup_{y \in K} \,\langle \lambda, y \rangle\right]\\[12pt]
&=&-\displaystyle\inf_{x\in \Re^n} {\cal L}(x,\lambda)+\delta_K^*(\lambda).
\end{array}
\end{equation}
Combining (\ref{eq:aa}) and (\ref{eq:bb}), we obtain $\nu^*(\lambda)=\theta(\lambda)$. From this, we can easily get the equality
$\nu^{**}(s)=\theta^*(s)$. This proves property (a).

The results in property (b) come from Theorem 11.1 of \cite{RW98}.

Now we prove property (c). Assume that $\nu$ is a proper lower semi-continuous function. For  $s \in {\rm dom}\, \nu$ with $\nu(s) \in \Re$, it follows from Theorem 11.1 of \cite{RW98} that $\nu^{**}(s)=\nu(s)$, implying $\nu (s)=\theta^*(s)$ from the equality
$\nu^{**}(s)=\theta^*(s)$ just proved. When $s \notin {\rm dom}\, \nu$ or $\nu(s)=+\infty$. Then $s \notin {\cal S}$. Let
$\hat s=\Pi_{{\cal S}}(s)$. Then $u=\hat s-s \ne 0$ and
$$
\forall (x,y)\in \Re^n \times K,\quad \langle y-g(x), u\rangle\geq \langle \hat s, u \rangle.
$$
Therefore, for any $\lambda \in {\rm dom}\,\theta$ and $t\geq 0$,
$$
\begin{array}{rcl}
\theta (\lambda-tu) &=&-\displaystyle\inf_{x\in \Re^n, y \in K} [f(x)+\langle \lambda-t\, u, g(x)-y \rangle]\\[8pt]
&\leq& \theta (\lambda)-t\, \langle \hat s, u\rangle.
\end{array}
$$
Then
$$
\begin{array}{rcl}
\theta^*(s) &=&\displaystyle{\sup_{\lambda'\in {\cal Y}}}\left\{ \langle \lambda', s\rangle-\theta (\lambda')\right\}\\[6pt]
&\geq& \langle s, \lambda-u\rangle-\theta (\lambda-tu)\\[6pt]
&\geq& \langle \lambda, s\rangle-\theta(\lambda)+t\,\|u\|^2.
\end{array}
$$
Since $t\geq 0$ is arbitrary and $u \ne 0$, $\theta^*(s)=+\infty$. Combining these two cases, we obtain $\nu(s)=\theta^*(s)$.
\hfill $\Box$
\begin{proposition}\label{prop:2.9}
 Let $s \in {\cal Y}$ and $\lambda \in {\cal Y}$. Assume that the function $\nu$ is proper lower semi-continuous with $\nu(s)\in \Re$. Then the following properties are equivalent:
\begin{itemize}
\item[{\rm (i)}] $s \in \partial \theta (\lambda)$;
\item[{\rm (ii)}] $\lambda \in \partial \nu^{**}(s)$ (actually $\lambda \in \partial \nu(s)$);
\item[{\rm (iii)}] $s \in {\cal S}$ and any solution to Problem (\ref{eq:p3y}) minimizes $l(\cdot,\cdot,\lambda)$ over $\Re^n \times {\cal Y}$;
    \item[{\rm (iv)}] There is a feasible  solution for Problem (\ref{eq:p3y}) that minimizes $l(\cdot,\cdot,\lambda)$ over $\Re^n \times {\cal Y}$.
\end{itemize}
\end{proposition}
{\bf Proof}. Since $\nu$ is proper lower semi-continuous convex, we have from Theorem 11.1 of \cite{RW98} and $\nu^*=\theta$
that $\theta$ is proper lower semi-continuous convex. Thus, we have from Theorem 23.5 of \cite{Rock70} that
\begin{equation}\label{eq:2.19}
s \in \partial \theta (\lambda) \Longleftrightarrow \theta (\lambda)+\theta^*(s)=\langle \lambda, s\rangle.
\end{equation}
[(i)$\Longleftrightarrow$(ii)] It follows from Corollary 23.5.1 of \cite{Rock70}
\footnote{Corollary 23.5.1. If $f$ is a closed proper convex function, $\partial f^*$ is the
inverse of $\partial f$ in the sense of multivalued mappings, i.e. $x \in \partial f^*(x^*)$ if and
only if  $x^*\in \partial f(x)$.} that
$$
s \in \partial \theta (\lambda) \Longleftrightarrow \lambda \in \partial \theta^*(s).
$$
By Proposition \ref{prop:2.8},
$$
s \in \partial \theta (\lambda) \Longleftrightarrow s \in \partial \nu^*(\lambda)\Longleftrightarrow
\lambda \in \partial \nu^{**}(s).
$$
Since $\nu$ is lower semi-continuous at $s$ with $\nu(s) \in \Re$, one has $\partial \nu(s)=\partial \nu^{**}(s)$. Thus we obtain
$s \in \partial \theta (\lambda)$ if and only if $\lambda \in \partial \nu (s)$.

[(i), (ii)$\Longrightarrow$(iii)] Let $s \in \partial \theta(\lambda)$. By (ii), $s \in {\rm dom}\,\nu^{**}={\rm dom}\,\nu={\cal S}$. Let $(x_s,y_s)$ be an arbitrary solution to Problem (\ref{eq:p3y}). Then
$$
\begin{array}{rcl}
l(x_s,y_s,\lambda)&=&f(x_s)-\langle \lambda, s\rangle \quad (g(x_s)+s=y_s)\\[6pt]
&=&\nu (s)-\langle \lambda, s\rangle \quad (\mbox{definition of } \nu)\\[6pt]
&=&\theta^*(s)-\langle \lambda, s\rangle \quad (\mbox{Proposition} \ref{prop:2.8})\\[6pt]
&=&-\theta (\lambda) \quad ((\ref{eq:2.19}) \mbox{ and } s \in \partial \theta (\lambda))\\[6pt]
&=&-\displaystyle\inf_{x\in \Re^n, y \in K} l(x,y,\lambda) \quad (\mbox{definition of } \theta).
\end{array}
$$
This implies that $(x_s,y_s)$ is a solution to the problem minimizing $l(\cdot,\cdot,\lambda)$ on $\Re^n \times K$.

[(iii)$\Longrightarrow$(iv)] This implication comes from the fact that Problem (\ref{eq:p3y}) has a solution when $s \in {\cal S}$ and ${\rm dom}\, \theta \ne \emptyset$ (this is from the fact that $\theta$ is a proper lower semi-continuous convex function).

[(iv)$\Longrightarrow$(i)] Let $(x_s,y_s)$ be a feasible point of  Problem (\ref{eq:p3y}) that minimizes $l(\cdot,\cdot,\lambda)$ on $\Re^n \times K$. For any $u \in {\cal Y}$,
$$
\begin{array}{rcl}
\langle s, u\rangle-\theta(u)&\leq& \langle s, u\rangle+ f(x_s)+\langle u, g(x_s)-y_s\rangle\quad (\mbox{definition of } \theta)\\[4pt]
&=&f(x_s)+\langle u, g(x_s)-y_s+s\rangle\\[4pt]
&=&f(x_s)+\langle \lambda, g(x_s)-y_s+s\rangle \quad (\mbox{feasibility of } (x_s,y_s), \mbox{ implying } g(x_s)+s=y_s)\\[4pt]
&=&\langle s, \lambda\rangle+\displaystyle \inf_{x\in \Re^n, y \in K} \left[f(x)+\langle \lambda, g(x)-y\rangle\right]\\[4pt]
&=&\langle s, \lambda\rangle-\theta (\lambda),
\end{array}
$$
which gives
$$
\theta (u)\geq \theta(\lambda)+\langle s, u-\lambda\rangle.
$$
Thus we obtain that $s \in\partial \theta (\lambda)$.\hfill $\Box$
\begin{remark}\label{remark-subdiff}
 If $\nu$ is a proper convex function, we have from Theorem 11.1 of \cite{RW98} and $\nu^*=\theta$
that $\theta$ is proper lower semi-continuous convex.
In this case, ${\rm dom}\, \theta\ne \emptyset$, we have from
 Corollary 23.5.1 of \cite{Rock70} that
 $$
 {\rm range}\, \partial \theta={\rm dom}\, \partial \theta^*.
 $$
 Moreover, if  $\nu$ is proper lower semi-continuous, then from Proposition \ref{prop:2.8}, we have
 $\nu=\theta^*$. Therefore, we obtain
$$
{\rm range}\, \partial \theta={\rm dom}\, \partial \theta^*={\rm dom}\, \partial \nu.
$$
It follows from Theorem 23.4 of \cite{Rock70} that
$$
{\rm ri}\,{\rm dom}\,\nu \subset {\rm dom}\, \partial \nu \subset {\rm dom}\,\nu.
$$
Therefore we obtain from ${\rm dom}\,\nu={\cal S}$ that
\begin{equation}\label{eq:Sri}
{\rm ri}\,{\cal S}\subset {\rm range}\, \partial \theta\subset {\cal S}.
\end{equation}
\end{remark}
\begin{proposition}\label{prop:p2.15}
Assume that $\nu$ is a proper lower semi-continuous function \rev{ and}
\begin{equation}\label{eq:domainE}
{\rm dom}\, \partial \nu=  {\rm dom}\,\nu.
\end{equation}
 Then ${\rm range}\, \partial \theta={\cal S}$.
\end{proposition}
{\bf Proof}. This result comes from Remark \ref{remark-subdiff} directly. \hfill $\Box$

The dual function associated with Problem ${\rm P}(s)$, denoted by $\theta_s:{\cal Y} \rightarrow
\overline \Re$, at $\lambda \in {\cal Y}$, is the value
$$
\begin{array}{rcl}
\theta_s(\lambda)&=&-\displaystyle \inf_{x\in \Re^n, y+s \in K} \left[f(x)+\langle \lambda, g(x)-y\rangle\right]\\[12pt]
&=&-\displaystyle \inf_{x\in \Re^n, y'\in K} \left[f(x)+\langle \lambda, g(x)+s-y'\rangle\right]\\[12pt]
&=&\theta (\lambda)-\langle s, \lambda\rangle.
\end{array}
$$
For an optimization problem ${\cal P}$, we use ${\rm val}\,({\cal P})$ and ${\rm Sol}\,({\cal P})$ to represent the optimal value and the set of optimal solutions of Problem ${\cal P}$. We use ${\rm D}$ and ${\rm D}(s)$ to denote the conjugate dual problems of Problems ${\rm P}$ and ${\rm P}(s)$, respectively. Then Problems ${\rm D}$ and ${\rm D}(s)$ can be expressed as follows:
\begin{equation}\label{p:d}
({\rm D}) \quad \max_{\lambda}\, [-\theta (\lambda)] \quad \quad \quad \quad \quad \quad ({\rm D}(s)) \quad \max_{\lambda}\, [\langle s, \lambda\rangle-\theta (\lambda)].
\end{equation}
\begin{proposition}\label{prop:facts}
We have
\begin{itemize}
\item[{\rm (i)}]${\rm val}\,({\rm P})=\nu(0)$, ${\rm val}\,({\rm D})=\nu^{**}(0)$, ${\rm val}\,({\rm D})=\theta^{*}(0)$, ${\rm Sol}\,({\rm D})=\partial \nu^{**}(0)$;
\item[{\rm (ii)}]If $\nu$ is lower semicontinuous at $0\in {\cal Y}$ with $\nu (0)$ being finite (in this case, Problem ${\rm P}$ is feasible), then
${\rm val}\,({\rm P})=\nu(0)=\nu^{**}(0)={\rm val}\,({\rm D})$, ${\rm Sol}\,({\rm D})=\partial \nu(0)$;
\item[{\rm (iii)}]For $s \in {\cal S}$, ${\rm val}\,({\rm P}(s))=\nu(s)$, ${\rm val}\,({\rm D}(s))=\nu^{**}(s)$, ${\rm val}\,({\rm D}(s))=\theta^{*}(s)$, ${\rm Sol}\,({\rm D}(s))=\partial \nu^{**}(s)$;
\item[{\rm (iv)}]For $s \in {\cal S}$, if $\nu$ is lower semicontinuous at $s$ with $\nu (s)$ being finite (in this case, Problem ${\rm P(s)}$ is feasible), then
${\rm val}\,({\rm P}(s))=\nu(s)=\nu^{**}(s)={\rm val}\,({\rm D}(s))$,  ${\rm Sol}\,({\rm D}(s))=\partial \nu(s)$;
\item[{\rm (v)}]If $s \in {\rm ri}\,{\cal S}$ with $\nu (s)$ being finite, then
${\rm val}\,({\rm P}(s))=\nu(s)=\nu^{**}(s)={\rm val}\,({\rm D}(s))$,  ${\rm Sol}\,({\rm D}(s))=\partial \nu(s)\ne \emptyset$;
\item[{\rm (vi)}]If $s \in {\rm int}\,{\cal S}$, then ${\rm val}\,({\rm P}(s))$ is finite and
${\rm val}\,({\rm P}(s))=\nu(s)=\nu^{**}(s)={\rm val}\,({\rm D}(s))$,  and ${\rm Sol}\,({\rm D}(s))=\partial \nu(s)$ is an nonempty compact set.
\end{itemize}
\end{proposition}
{\bf Proof}. (i)-(iv) come from Proposition 2.118 of \cite{BS00}. We only need to prove (v) and (vi).
For $s \in {\rm ri}\,{\cal S}$, there exists $x^0$ such that $g(x^0)+s \in {\rm ri}\, K$, which implies that the generalized Slater condition holds for Problem ${\rm P}(s)$:
$$
  \begin{array}{cl}
  \min & f(x)\\
{\rm s.t.}& g(x)+s\in K.
\end{array}
$$
If follows from (2.311) of \cite{BS00} that
$$
{\rm dom} \, \nu = K -g({\rm dom}\,f).
$$
The generalized Slater condition implies that
$$
s \in {\rm ri}\, {\rm dom}\,\nu.
$$
Therefore, we have from Theorem 23.4 of \cite{Rock70} that $\partial \nu (s)\ne \emptyset$. It follows from Proposition 2.118 of \cite{BS00} that $\nu(s)=\nu^{**}(s)$ and $\partial \nu^{**}(s)=\partial \nu(s)$.  Thus we obtain all results in (v).

If $s \in {\rm int}\,{\cal S}$, from the above analysis, we obtain
$$
s \in {\rm int}\, {\rm dom}\,\nu.
$$
The results in (vi) can be obtained by using Theorem 23.4 of \cite{Rock70}. The proof is completed. \hfill $\Box$
\section{The augmented Lagrangian method}\label{sec:ALM}
\setcounter{equation}{0}
In this section, we focus on the convex optimization problem with the least constraint violation for the case that $0 \notin {\cal S}$, \rev{where ${\cal S}$ is a closed set}; namely, the following problem
 \begin{equation}\label{eq:1.12a}
 {\rm P}(\bar s)\quad \quad
 \begin{array}{ll}
  \min & f(x)\\
{\rm s.t.}& g(x)+\bar s\in K,
\end{array}
\end{equation}
where $\bar s$ is the projection of $0 \in {\cal Y}$ on ${\cal S}$; namely,
  \begin{equation}\label{eq:snshifta}
  \bar s=\mbox{argmin} \left\{\displaystyle \frac{1}{2}\|s\|^2: s \in {\cal S}\right\}.
  \end{equation}
  If ${\cal S}$ is closed and $0 \notin {\cal S}$, then $\bar s$ is on the relative boundary of ${\cal S}$; namely, $\bar s\in {\rm ribdry}\,{\cal S}$. Now we present conditions ensuring the zero duality gap for ${\rm P}({\bar s})$ and its conjugate dual
  ${\rm D}({\bar s})$, as well as the characterization of ${\rm Sol}\,{\rm D}({\bar s})$.
  \begin{proposition}\label{prop:bars}
  Assume that $\bar s \ne 0$ and ${\rm val}\,{\rm P}(\bar s)\in \Re$.
  \begin{itemize}
  \item[{\rm (i)}]
  Suppose that $\nu$ is lower semi-continuous at $\bar s$. Then
  $$
  {\rm val}\,{\rm D}(\bar s)={\rm val}\,{\rm P}(\bar s).
  $$
  \item[{\rm (ii)}]If $\partial \nu (\bar s)\ne \emptyset$, then
  $$
  {\rm val}\,{\rm D}(\bar s)={\rm val}\,{\rm P}(\bar s)\quad \mbox{ and } \quad {\rm Sol}\,{\rm D}(\bar s)=\partial \nu (\bar s)
  $$
  or
  $$
  {\rm Sol}\,{\rm D}(\bar s)=\{\lambda \in {\cal Y}: \bar s \in \partial \theta (\lambda)\}=[\partial \theta]^{-1}(\bar s).
  $$
  \end{itemize}
  \end{proposition}
  \begin{proposition}\label{prop:unb}
  Assume that $\bar s \ne 0$, ${\rm val}\,{\rm P}(\bar s)\in \Re$,
 $\nu$ is lower semi-continuous at $\bar s$  and ${\rm Sol}\,{\rm D}(\bar s)\ne \emptyset$.
  Then ${\rm Sol}\,{\rm D}(\bar s)$ is unbounded with
  \begin{equation}\label{eq:recD}
  -\bar s\in [{\rm Sol}\,{\rm D}(\bar s)]^{\infty}.
  \end{equation}
  \end{proposition}
  {\bf Proof}. From Proposition \ref{prop:bars}, we have
  $$
  {\rm Sol}\,{\rm D}(\bar s)=\{\lambda \in {\cal Y}: 0\in \partial \theta_{\bar s}(\lambda)\}=\{\lambda \in {\cal Y}:\bar s\in \partial \theta (\lambda)\}.
  $$
  Then, for any $\bar \lambda \in {\rm Sol}\,{\rm D}(\bar s)$, we have that
  $$
  [{\rm Sol}\,{\rm D}(\bar s)]^{\infty} =\{\xi \in {\cal Y}: 0\in \partial \theta_{\bar s}(\bar \lambda+t\xi),\forall t\geq 0\}
  =\{\xi \in {\cal Y}: \bar s \in \partial \theta (\bar \lambda+t\xi),\forall t\geq 0\}.
  $$
  Therefore we only need to prove
  \begin{equation}\label{eq:chk}
  \bar s \in \partial \theta (\bar \lambda-t \bar s), \quad  \forall t\geq 0.
  \end{equation}
  From the definition of $\bar s$, one has that
  $$
  \bar s=\mbox{argmin} \left\{\displaystyle \frac{1}{2}\|s\|^2: s \in {\cal S}\right\},
 $$
  where ${\cal S}$ is a convex set in ${\cal Y}$ defined by (\ref{eq:n1}).  Then one has
  $$
  \langle \bar s, s-\bar s\rangle\geq 0\rev{,} \quad \forall s \in {\cal S}
  $$
  or
  \begin{equation}\label{eq:10000}
  \langle \bar s, s\rangle\geq \|\bar s\|^2\rev{,} \quad \forall s \in {\cal S}.
  \end{equation}
  For any $\lambda \in {\rm dom}\, \theta$, we have for any $t \geq 0$ that
  $$
  \begin{array}{rcl}
  \theta (\lambda)-\theta(\bar \lambda-t\bar s)
  &=&\theta (\lambda)+\displaystyle \inf_{x\in \Re^n, y\in K} \left[f(x)+\langle \bar \lambda-t\bar s, g(x)-y\rangle\right]\\[10pt]
  &\geq& \theta (\lambda)+\displaystyle \inf_{x\in \Re^n, y\in K} \left[f(x)+\langle \bar \lambda, g(x)-y\rangle\right]+\displaystyle \inf_{x\in \Re^n, y\in K} \left[f(x)+\langle -t\bar s, g(x)-y\rangle\right]\\[10pt]
  &=&\theta (\lambda)-\theta (\bar \lambda)+\displaystyle \inf_{s \in {\cal S}}\left[t\,\langle \bar s,s\rangle\right]\\[10pt]
  &=&\theta (\lambda)-\theta (\bar \lambda)+t\displaystyle \inf_{s \in {\cal S}}\,\langle \bar s,s\rangle \\[10pt]
  &\geq& \langle \bar s, \lambda-\bar \lambda \rangle+t\, \langle \bar s,\bar s\rangle \quad \quad (\mbox{from }(\ref{eq:10000})\mbox{ and } \bar s\in \partial \theta (\bar \lambda))\\[8pt]
  &=&\langle \bar s, \lambda-(\bar \lambda -t\bar s)\rangle,
  \end{array}
  $$
  implying the truth of (\ref{eq:chk}).  The proof is completed. \hfill $\Box$

  The analysis of the augmented Lagrangian method is based on the notion of the Moreau-Yosida regularization and the proximal mapping of convex functions.  Let $\theta_r: {\cal Y} \rightarrow \overline \Re$ be the Moreau-Yosida regularization of $\theta$; namely,
\begin{equation}\label{eq:MY}
\theta_r (\lambda)=\inf_{\lambda'\in {\cal Y}} \left\{\theta (\lambda')+\displaystyle \frac{1}{2r}\|\lambda'-\lambda\|^2 \right\}.
\end{equation}
The proximal mapping $P_{r\theta}:{\cal Y} \rightarrow {\cal Y}$ is defined by
\begin{equation}\label{eq:PM}
P_{r\theta} (\lambda)=\mbox{argmin} \left\{\theta (\lambda')+\displaystyle \frac{1}{2r}\|\lambda'-\lambda\|^2: \lambda' \in {\cal Y} \right\}.
\end{equation}
\begin{lemma}\label{lem:2.4}
The dual function $\theta$ is lower semi-continuous convex. Suppose that $\theta$ is proper and let $r>0$. Then
\begin{equation}\label{eq:2.9}
\theta_r(\lambda)=-\inf_{x\in \Re^n, y \in K} l_r(x,y,\lambda),
\end{equation}
where $l_r$ is the augmented Lagrangian. Let
$
(x(\lambda,r), y(\lambda,r)) =\mbox{argmin}\left\{ l_r(x,y,\lambda): x\in \Re^n, y \in K \right\}.
$
Then
$$
P_{r\theta}(\lambda)=\lambda+r[g(x(\lambda,r))-y(\lambda,r)] \quad \mbox{ and } \quad y(\lambda,r)-g(x(\lambda,r))\in \partial
\theta (P_{r\theta}(\lambda)).
$$
\end{lemma}
\begin{proposition}\label{prop-2.5}
Let $\lambda \in Y$ and $r>0$. If $\lambda \in {\rm dom }\,\theta$, then the augmented Lagrangian subproblem
\begin{equation}\label{eq:2.9a}
\min_{x \in \Re^n, y \in K} l_r(x,y,\lambda)
\end{equation}
has a solution.
\end{proposition}
{\bf Proof}. Since $\lambda \in {\rm dom }\,\theta$, $ {\rm dom }\,\theta \ne \emptyset$, from (\ref{eq:2.7}), we have that $\theta$ is a proper lower semi-continuous convex function, and the optimal value $\theta_r(\lambda)$ of the problem in the right hand side  of (\ref{eq:MY}) is finite. By Lemma \ref{lem:2.4}, Problem (\ref{eq:2.9a}) has a solution.
\hfill $\Box$
  \begin{proposition}\label{prop:2.17}
   Assume that $\bar s \ne 0$, ${\rm val}\,{\rm P}(\bar s)\in \Re$,
 $\nu$ is lower semi-continuous at $\bar s$  and ${\rm Sol}\,{\rm D}(\bar s)\ne \emptyset$. Let $\lambda \in {\cal Y}$. Then the following properties hold:
 \begin{itemize}
 \item[${\rm (i)}$]${\rm dist}\, \left(\lambda-\alpha \bar s, {\rm Sol}\,{\rm D}(\bar s)\right)\leq {\rm dist}\, \left(\lambda, {\rm Sol}\,{\rm D}(\bar s)\right), \ \forall \alpha \geq 0$;
  \item[${\rm (ii)}$]$P_{r\theta}(\lambda)=P_{r\theta_{\bar s}}(\lambda-r \bar s)$;
   \item[${\rm (iii)}$]${\rm dist}\, (P_{r\theta}(\lambda),{\rm Sol}\,{\rm D}(\bar s))\leq {\rm dist}\, (\lambda, {\rm Sol}\,{\rm D}(\bar s))$.
 \end{itemize}
  \end{proposition}
  {\bf Proof}. (i). Let
  $$
  \tilde \lambda=\Pi_{{\rm Sol}\,{\rm D}(\bar s)}(\lambda),
  $$
  which is well-defined since ${\rm Sol}\,{\rm D}(\bar s)$ is an nonempty closed convex set. It follows from
  Proposition \ref{prop:unb} that $\tilde \lambda-\alpha \bar s \in {\rm Sol}\,{\rm D}(\bar s)$ for any $\alpha \geq 0$.
  Then
  $$
  {\rm dist}\, (\lambda-\alpha \bar s, {\rm Sol}\,{\rm D}(\bar s))\leq \|\lambda-\alpha \bar s-[\tilde \lambda-\alpha \bar s]\|=\|\lambda-\tilde \lambda\|={\rm dist}\, (\lambda, {\rm Sol}\,{\rm D}(\bar s)).
  $$
  (ii). Let
  $u=P_{r\theta_{\bar s}}(\lambda-r \bar s)$. Then
  $$
  0 \in \partial \theta_{\bar s}(u)+\displaystyle \frac{1}{r}[u-(x-r \bar s)].
  $$
  From this, there exists some $\tilde s\in \partial \theta_{\bar s}(u)$ such that
  $$
  0 =\tilde s+\displaystyle \frac{1}{r}[u-(x-r \bar s)]
  $$
  or equivalently
  \begin{equation}\label{eq:star}
  u=\lambda-r(\tilde s+\bar s).
  \end{equation}
  From the expression $\partial \theta_{\bar s}(\lambda)=\partial \theta (\lambda)-\bar s$, we get
  $$
  \tilde s+\bar s\in \partial \theta_{\bar s}(u)
  $$
  and from (\ref{eq:star}) that
  $$
  0 \in \partial \theta (u)+\displaystyle \frac{1}{r}[u-\lambda].
  $$
  Therefore we obtain $u =P_{r\theta}(\lambda)$.
  This proves (ii).\\
  (iii). In view of (ii), we have
  \begin{equation}\label{eq:2.33}
  {\rm dist}\, \left(P_{r\theta}(\lambda), {\rm Sol}\,{\rm D}(\bar s)\right)= {\rm dist}\,\left(P_{r\theta_{\bar s}}(\lambda-\alpha \bar s), {\rm Sol}\,{\rm D}(\bar s)\right).
  \end{equation}

  Since ${\rm Sol}\,{\rm D}(\bar s)={\rm argmin}\, \theta_{\bar s}$ and
  $$
  {\rm argmin}\, \theta_{\bar s}=\{u \in {\cal Y}: u=P_{r\theta_{\bar s}}(u)\},
  $$
  in view of the fact that $P_{r\theta_{\bar s}}$ is an nonexpansive mapping, we have for any $u \in {\rm Sol}\,{\rm D}(\bar s)$ that
  $$
  \|P_{r\theta_{\bar s}}(\lambda-\alpha \bar s)-u\|=\|P_{r\theta_{\bar s}}(\lambda-\alpha \bar s)-P_{r\theta_{\bar s}}(u)\|
  \leq \|(\lambda-\alpha \bar s)-u\|.
  $$
  This implies that a proximal step decreases the distance to the set of minimizers:
  \begin{equation}\label{eq:2.34}
   {\rm dist}\,\left(P_{r\theta_{\bar s}}(\lambda-\alpha \bar s), {\rm Sol}\,{\rm D}(\bar s)\right)
   \leq {\rm dist}\,\left(\lambda-\alpha \bar s,{\rm Sol}\,{\rm D}(\bar s)\right).
  \end{equation}
  The inequality in (iii)  is now achieved by combining (\ref{eq:2.33}), (\ref{eq:2.34}) and (i). \hfill $\Box$\\
  \begin{lemma}\label{lem:2.13}
  Assume that $\bar s \ne 0$. Assume also that $g:\Re^n \rightarrow {\cal Y}$ is a smooth mapping and $G_K$ is a graph-convex mapping. Then the following properties of $(\bar x,\bar y) \in \Re^n \times {\cal Y}$ are equivalent:
  \begin{itemize}
  \item[{\rm (i)}]$\bar y-g(\bar x)=\bar s$ and $\bar y \in K$;
    \item[{\rm (ii)}]${\rm D}g(\bar x)^*(g(\bar x)-\bar y)=0$ and $\Pi_K(g(\bar x))=\bar y$;
  \item[{\rm (iii)}]$(\bar x,\bar y)$ is a solution to
  \begin{equation}\label{eq:2.23}
  \displaystyle \min_{x\in \Re^n, y\in K} \displaystyle \frac{1}{2}\|g(x)-y\|^2.
  \end{equation}
  \end{itemize}
  \end{lemma}
  {\bf Proof}. It suffices to prove that (iii) is equivalent (i) and (iii) is equivalent (ii). By introducing
  \begin{equation}\label{eq:gy}
  s=-g(x)+y,
  \end{equation}
  Problem (\ref{eq:2.23}) is equivalent to
  \begin{equation}\label{eq:2.23a}
  \begin{array}{cl}
   \displaystyle \min_{x\in \Re^n, s\in {\cal Y}} & \displaystyle \frac{1}{2}\|s\|^2\\[8pt]
   {\rm s.t.} & g(x)+s \in K.
   \end{array}
  \end{equation}
  This implies the equivalence between (i) and (iii).
  Since $G_K$ is a graph-convex mapping, we have that Problem (\ref{eq:2.23a}) is a  convex optimization problem.
    Let $(\bar x, \bar s)$ be a solution to Problem (\ref{eq:2.23a}). Noting that the generalized Slater condition holds for  Problem (\ref{eq:2.23a}), we have that $(\bar x, \bar s)$ is a solution to Problem (\ref{eq:2.23a}) if and only if the following KKT conditions hold at $(\bar x, \bar s)$; namely, there exists a Lagrangian multiplier $\bar \lambda \in {\cal Y}$ such that
    \begin{equation}\label{psKKT}
    \begin{array}{l}
    {\rm D}g(\bar x)^*\bar \lambda=0,\\[4pt]
    \bar s+\bar \lambda=0,\\[4pt]
    \bar \lambda\in N_K(g(\bar x)+\bar s).
    \end{array}
    \end{equation}
  Letting $\bar y=g(\bar x)+\bar s$, the above relations are equivalently expressed as follows
  \begin{equation}\label{psKKTa}
    \begin{array}{l}
    {\rm D}g(\bar x)^*(g(\bar x)-\bar y)=0,\\[4pt]
        g(\bar x)-\bar y \in N_K(\bar y).
    \end{array}
    \end{equation}
  Noting that $g(\bar x)-\bar y \in N_K(\bar y)$ is equivalent to $\bar y=\Pi_K(g(\bar x))$, we obtain that (\ref{psKKTa}) is equivalent to (i). The proof is completed. \hfill $\Box$

  Now we provide a set of optimality conditions for characterizing a solution for the optimization with the least constraint violation ${\rm P}(\bar s)$ in terms of the augmented Lagrangian.
  \begin{theorem}\label{Optimality-AL}
  Assume that $\nu$ is a proper lower semi-continuous function  with $\nu(\bar s)\in \Re$ and
\begin{equation}\label{eq:domainEa}
\bar s \in {\rm dom}\, \partial \nu.
\end{equation}
  Let $r>0$, $l_r$ be the augmented Lagrangian defined by (\ref{eq:2.5agm}).  Assume that $g$ is a smooth mapping from $\Re^n$ to ${\cal Y}$. Then $(\bar x, \bar y)$ is a solution to the following problem
  \begin{equation}\label{eq:1.12as}
 \begin{array}{ll}
  \min & f(x)\\[3pt]
{\rm s.t.}& g(x)+\bar s=y,\\[3pt]
   & y \in K
\end{array}
\end{equation}
if and only if there exists some $\bar \lambda \in {\cal Y}$ such that
\begin{equation}\label{eq:23567}
\begin{array}{l}
(\bar x,\bar y) \in \displaystyle {\rm argmin}_{x\in \Re^n, y \in K} l_r(x,y,\bar \lambda),\\[8pt]
{\rm D}g(\bar x)^*(g(\bar x)-\bar y)=0,\\[6pt]
\Pi_K(g(\bar x))=\bar y.
\end{array}
\end{equation}
  \end{theorem}
  {\bf Proof}. Necessity. Let $(\bar x, \bar y)$ be a solution to Problem
  (\ref{eq:1.12as}). Then $\bar y-g(\bar x)=\bar s$, $\bar y\in K$, where
  $$
  \bar s=\mbox{argmin} \left\{\displaystyle \frac{1}{2}\|s\|^2: s \in {\cal S}\right\}.
  $$
  Then by the implication (i)$\Longrightarrow$ (ii) of Lemma \ref{lem:2.13}, we obtain
  $$
  {\rm D}g(\bar x)^*(g(\bar x)-\bar y)=0;\,\,
\Pi_K(g(\bar x))=\bar y,
  $$
  namely, the second and third relations in (\ref{eq:23567}) are valid.

    It follows from the equivalence of (i) and (ii) in Proposition \ref{prop:2.9} that $\bar s \in  {\rm dom}\, \partial \nu$
    implies that there exists some $\bar \lambda$ such that $\bar s \in \partial \theta (\bar \lambda)$. By the implication
    (i)$\Longrightarrow$ (iii) of Proposition \ref{prop:2.9}, $(\bar x,\bar y)$ minimizes the Lagrangian $l(\cdot,\cdot,\bar \lambda)$ over $\Re^n \times K$:
    \begin{equation}\label{eq:2.38}
    f(\bar x)+\langle \bar \lambda, g(\bar x)-\bar y\rangle\leq f(x)+\langle \bar \lambda, g(x)-y\rangle\rev{,}\quad \forall (x,y) \in \Re^n \times K.
    \end{equation}
    For any $(x,y) \in \Re^n \times K$, $y-g(x) \in {\cal S}$ so that
    $$
    \|g(\bar x)-\bar y\|=\|\bar s\|\leq \|g(x)-y\|
    $$
    by the definition of $\bar s$. Using (\ref{eq:2.38}), we get for any $(x,y) \in \Re^n \times K$ that
    $$
    f(\bar x)+\langle \bar \lambda, g(\bar x)-\bar y\rangle+\displaystyle
    \frac{r}{2}\|g(\bar x)-\bar y\|^2
    \leq f(x)+\langle \bar \lambda, g(x)- y\rangle+\displaystyle
    \frac{r}{2}\|g(x)-y\|^2,
    $$
    which proves
    $$
    (\bar x,\bar y) \in \displaystyle {\rm argmin}\,\Big \{l_r(x,y,\bar \lambda):\,x\in \Re^n, y \in K\Big\}
    $$
    in (\ref{eq:23567}).

    Sufficiency. In view of the implication (ii) $\Longrightarrow$ (i) of Lemma \ref{lem:2.13}, ${\rm D}g(\bar x)^*(g(\bar x)-\bar y)=0$ and $\Pi_K(g(\bar x))=\bar y$ imply that $(\bar x,\bar y)$ satisfies the constraints of Problem (\ref{eq:1.12as}). Now let $(x,y)$ satisfy $g(x)=y, y+\bar s \in K$. Then by
    $$
     (\bar x,\bar y) \in \displaystyle {\rm argmin}\,\Big\{ l_r(x,y,\bar \lambda):\,x\in \Re^n, y \in K\Big\}
    $$
    and $g(\bar x)-\bar y=g(x)-y=-\bar s$, we have that
    $$
    f(\bar x)-\langle \bar \lambda, \bar s\rangle+\displaystyle \frac{r}{2}\|\bar s\|^2 \leq f(x)-\langle \bar \lambda, \bar s\rangle+\displaystyle \frac{r}{2}\|\bar s\|^2.
    $$
    Hence we get $f(\bar x)\leq f(x)$ for all $(x,y)$ satisfying $g(x)=y$ and $y+\bar s\in K$, which implies $(\bar x,\bar y)$ is a solution to Problem (\ref{eq:1.12as}). The proof is completed. \hfill $\Box$

    Now we are ready to describe the augmented Lagrangian method for Problem (\ref{eq:p2}).

    \begin{algorithm}[h]
	\caption{The Augmented Lagrangian Method} \label{al1}
	Initialize $\lambda^{0}\in{\cal Y}$ and $r_{0}>0$. Set $k:=0$.\\
	\While{the stopping condition does not hold}{
\begin{itemize}
\item[1.] Find a solution of
$$
\min_{x\in \Re^n, y \in K} l_{r_k}(x,y,\lambda^k)
$$
and denote it by $(x^{k+1},y^{k+1})$.
\item[2.] Update the multiplier
$$
\lambda^{k+1}=\lambda^k+r_k(g(x^{k+1})-y^{k+1}).
$$
\item[3.] Choose a new penalty parameter $r_{k+1}\geq r_k$.
\end{itemize}
				Set $k:=k+1$.
	}
\end{algorithm}
\mbox{}\\[-1pt]
Define
\begin{equation}\label{eq:skD}
s^k=y^k-g(x^k).
\end{equation}
Then, from Lemma \ref{lem:2.4}, we have
\begin{equation}\label{eq:rls}
s^{k+1}\in \partial \theta (\lambda^{k+1}) \quad \mbox{ and }\quad \lambda^{k+1}=P_{r_k\theta}(\lambda^k).
\end{equation}
\begin{theorem}\label{prop:3.1}
Assume that $\bar s \ne 0$, ${\rm val}\,{\rm P}(\bar s)\in \Re$,
 $\nu$ is lower semi-continuous at $\bar s$  and ${\rm Sol}\,{\rm D}(\bar s)\ne \emptyset$.
 Let $\{(x^k,y^k,\lambda^k)\}$ be generalized by the augmented Lagrangian method.
 Then
 \begin{itemize}
 \item[{\rm (i)}] The sequence $\{\|s^k\|\}$ is nonincreasing;
 \item[{\rm (ii)}] The sequence $\{{\rm dist}\, (\lambda^k, {\rm Sol}\,{\rm D}(\bar s))\}$ is nonincreasing;
 \item[{\rm (iii)}] If $r_k \geq \underline{r}$ for some $\underline{r}>0$, then $s^k \rightarrow \bar s$.
 \end{itemize}
\end{theorem}
{\bf Proof}. (i). Noting that $s^k \in \partial \theta (\lambda^k)$,
$s^{k+1} \in \partial \theta (\lambda^{k+1})$ and
$$s^{k+1}=y^{k+1}-g(x^{k+1})=\displaystyle \frac{1}{r_k}[\lambda^k-\lambda^{k+1}],
$$
 we obtain
 $$
 \begin{array}{rcl}
 \|s^k\|^2&=&\|(s^k-s^{k+1})+s^{k+1}\|^2\\[6pt]
 &=&\|s^k-s^{k+1}\|^2+2\langle s^k-s^{k+1},s^{k+1}\rangle+\|s^{k+1}\|^2\\[6pt]
 &=&\|s^k-s^{k+1}\|^2+\displaystyle \frac{2}{r_k}\langle s^k-s^{k+1},\lambda^k-\lambda^{k+1}\rangle+\|s^{k+1}\|^2\\[6pt]
 &\geq& \|s^{k+1}\|^2,
 \end{array}
 $$
 in which
 $
 \displaystyle \frac{2}{r_k}\langle s^k-s^{k+1},\lambda^k-\lambda^{k+1}\rangle \geq 0
 $
 is used. This proves (i).

 (ii). In view of Lemma \ref{lem:2.4}, we have $\lambda^{k+1}=P_{r_k\theta}(\lambda^k)$. Applying Proposition \ref{prop:2.17}, we obtain (ii).

 (iii). It follows from Proposition \ref{prop:unb} that, if ${\rm Sol}\,{\rm D}(\bar s)\ne \emptyset$, then $-\bar s \in {\rm Sol}\,{\rm D}(\bar s)^{\infty}$.  Define a sequence in ${\rm Sol}\,{\rm D}(\bar s)$ as follows:
 $$
 u^0 \in {\rm Sol}\,{\rm D}(\bar s)\ \mbox{ and }\ u^{k+1}=u^k-r_k\bar s \quad (\forall k \geq 0).
 $$
Then we have $\{u^k\} \subset {\rm Sol}\,{\rm D}(\bar s)$.  Since $\lambda^{k+1}=\lambda^k-r_k s^{k+1}$ (from the augmented Lagrangian method), one has that
\begin{equation}\label{eq:ccc}
\lambda^k-u^k=\lambda^{k+1}-u^{k+1}+r_k(s^{k+1}-\bar s).
\end{equation}
Noting that $s^{k+1} \in \partial \theta (\lambda^{k+1})$ and $\bar s \in \partial \theta (u^{k+1})$ by Proposition \ref{prop:bars} (ii), the monotonicity of $\partial \theta$  implies that
$$
\langle s^{k+1}-\bar s,\lambda^{k+1}-u^{k+1}\rangle\geq 0.
$$
Thus, taking the square norm of the both sides of (\ref{eq:ccc}) and neglecting the term $\langle s^{k+1}-\bar s,\lambda^{k+1}-u^{k+1}\rangle$ in the right hand side yield
\begin{equation}\label{eq:ddd}
\|\lambda^k-u^k\|^2 \geq \|\lambda^{k+1}-u^{k+1}\|^2+r_k^2\|s^{k+1}-\bar s\|^2.
\end{equation}
This implies that the nonnegative sequence $\{\|\lambda^k-u^k\|\}$ is nonincreasing and hence converges. Therefore, we obtain from (\ref{eq:ddd}) that $r_k^2\|s^{k+1}-\bar s\|^2$ converges to zero. Since $r_k \geq \underline{r}$ for some $\underline{r}>0$, we have that $s^k \rightarrow \bar s$. \hfill $\Box$
\begin{corollary}\label{cor:unb-multiplier}
Assume that $\bar s \ne 0$, ${\rm val}\,{\rm P}(\bar s)\in \Re$,
 $\nu$ is lower semi-continuous at $\bar s$  and ${\rm Sol}\,{\rm D}(\bar s)\ne \emptyset$.
 Let $\{(x^k,y^k,\lambda^k)\}$ be generalized by the augmented Lagrangian method with $r_k \geq \underline r$ for some $\underline r>0$.
 Then $\{\lambda^k\}$  diverges.
\end{corollary}
{\bf Proof.} Noting
$$
\|\lambda^{k+1}-\lambda^k\|=\|-r_ks^{k+1}\|\geq \underline{r}\, \|s^{k+1}\| \rightarrow \underline{r}\,\|\bar s\|>0,
$$
we obtain that $\{\lambda^k\}$ is divergent. \hfill $\Box$

Define the set of accumulation points of $\{(x^k,y^k)\}$:\footnote{Here $\limsup$ stands for the outer limit of a sequence of sets from Chapter 4 of \cite{RW98}:
$$
\limsup_{k \rightarrow +\infty} C^k=\Big\{z: \mbox{there exists a subsequence } N\subset \textbf{N}, \exists\, z^k \in C^k
\mbox{ for } k \in N \mbox{ such that } z^k\stackrel{N} \rightarrow z\Big\}.
$$
}
$$
\omega=\displaystyle \limsup_{k \rightarrow \infty} \{(x^k,y^k)\}.
$$
\begin{proposition}\label{prop:sc}
Assume that $\bar s \ne 0$, ${\rm val}\,{\rm P}(\bar s)\in \Re$,
 $\nu$ is lower semi-continuous at $\bar s$  and ${\rm Sol}\,{\rm D}(\bar s)\ne \emptyset$.
Suppose that $\omega \ne \emptyset$. Then for any $(\bar x,\bar y)\in \omega$,
\begin{equation}\label{eq:4.5}
{\rm D}g(\bar x)^*(g(x^k)-y^k) \rightarrow 0\quad \mbox{ and }\quad \Pi_K(g(x^k))-y^k \rightarrow 0.
\end{equation}
\end{proposition}
{\bf Proof}. For any $(\bar x,\bar y) \in \omega$,  there exists a subsequence $ N\subset \textbf{N}$ such that
$(x^k,y^k)\stackrel{N} \rightarrow (\bar x,\bar y)$.  Since $r_k \geq \underline{r}$, by Theorem \ref{prop:3.1}, $s^k=y^k-g(x^k) \rightarrow \bar s$, we have $\bar y-g(\bar x)=\bar s$, $\bar y \in K$ and
$$
  \bar s=\mbox{argmin} \left\{\displaystyle \frac{1}{2}\|s\|^2: s \in {\cal S}\right\}.
  $$
It follows from the equivalence of (i) and (ii) of Lemma \ref{lem:2.13} that
\begin{equation}\label{eq:staraa}
{\rm D}g(\bar x)^* \bar s=0, \ \ \Pi_K(g(\bar x))-\bar y=0.
\end{equation}
From (\ref{eq:staraa}), we can easily obtain ${\rm D}g(\bar x)^*(g(x^k)-y^k) \rightarrow 0$; namely, the first property in (\ref{eq:4.5}) holds.

Let us denote
$$
\tilde y^k =\Pi_K(g(x^k)).
$$
Then we have from properties of the projection that
$$
\langle \tilde y^k-g(x^k),y-\tilde y^k \rangle\geq 0\rev{,} \quad \forall y \in K.
$$
Taking $y=y^k \in K$ in the above inequality yields
\begin{equation}\label{eq:4.6}
\langle \tilde y^k-g(x^k),y^k-\tilde y^k \rangle\geq 0.
\end{equation}
Notice that $\bar s =\Pi_{{\cal S}}(0)$ is characterized by
$$
\langle \bar s, s -\bar s \rangle\geq 0\rev{,} \quad \forall s \in {\cal S}.
$$
Taking $s=\tilde y^k-y^k+s^k=\tilde y^k-g(x^k) \in {\cal S}$, we obtain
\begin{equation}\label{eq:4.7}
\langle \bar s, \tilde y^k-y^k+s^k -\bar s \rangle\geq 0.
\end{equation}
Adding (\ref{eq:4.6}) and (\ref{eq:4.7}), we get
$$
\langle \bar s- \tilde y^k+g(x^k), \tilde y^k-y^k\rangle+\langle \bar s, s^k -\bar s \rangle\geq 0.
$$
Using $s^k=y^k-g(x^k)$ and the Cauchy-Schwartz inequality,
$$
\begin{array}{rcl}
\|\tilde y^k-y^k\|^2 & \leq& \|\tilde y^k-y^k\|^2+\langle \bar s- \tilde y^k+y^k-s^k, \tilde y^k-y^k\rangle+\langle \bar s, s^k -\bar s \rangle\\[6pt]
&=&\langle \bar s-s^k, \tilde y^k-y^k\rangle+\langle \bar s, s^k -\bar s \rangle\\[6pt]
&\leq& \|\bar s-s^k\|\,\| \tilde y^k-y^k\|+\|\bar s\|\,\|s^k -\bar s\|.
\end{array}
$$
Since $s^k \rightarrow \bar s$, the above inequality implies that there exists a constant $\beta>0$ such that
$$
\|\tilde y^k-y^k\|\leq \beta\, \|s^k-\bar s\|^{1/2}
$$
for sufficiently large $k$, which implies the  property $\Pi_K(g(x^k))-y^k \rightarrow 0$.  The proof is completed. \hfill $\Box$\\

Now we are in a position to state the main result in this paper, which shows that the augmented Lagrangian method can find an approximate solution to the optimization problem with the least constraint violation.
\begin{theorem}\label{epsilon-solution}
Consider the augmented Lagrangian method for Problem (\ref{eq:p2}).
Assume that $\bar s \ne 0$, ${\rm val}\,{\rm P}(\bar s)\in \Re$,
 $\nu$ is lower semi-continuous at $\bar s$  and ${\rm Sol}\,{\rm D}(\bar s)\ne \emptyset$. Suppose that $\omega \ne \emptyset$. Then for every $\varepsilon>0$, there exists a subsequence $ N\subset \textbf{N}$ such that
 \begin{equation}\label{eq:approx-opt1}
  (x^k,y^k) \in {\rm argmin}\, l_{r_{k-1}}(x,y,\lambda^{k-1})
  \end{equation}
  and
 \begin{equation}\label{eq:approx-opt2}
 \|{\rm D}g(x^k)^*(g(x^k)-y^k)\| \leq \varepsilon \quad \mbox{ and }\quad
  \|\Pi_K(g(x^k))-y^k\|\leq \varepsilon
  \end{equation}
 for every $k \in N$.
\end{theorem}
{\bf Proof}. The property (\ref{eq:approx-opt1}) comes from the definition of $(x^k,y^k)$ in the augmented Lagrangian method. The property (\ref{eq:approx-opt2}) follows  from Proposition \ref{prop:sc}.  The proof is completed. \hfill $\Box$
\section{Discussions}
\setcounter{equation}{0}
There are many practical backgrounds for the importance dealing with optimization problems with least constraints violations when the feasible sets of problems are possibly empty. Dai and Zhang \cite{DZhang2021} established necessary optimality conditions by reformulating them as optimization problems with an Lipschitz continuous equality or optimization problems with complementarity constraints. They also analyzed a penalty function method and a smoothing function method for solving these  optimization problems with the least constraint violation.

Can we construct an algorithm for solving a constrained optimization with the property that the method finds a solution when the problem is feasible, and finds a solution to the problem with the least constraint violation when the problem is infeasible? The work by  Chiche and  Gilbert\cite{ChCh2016} gave us a positive answer to the question for the convex quadratic optimization and they found that the augmented Lagrangian method is such a method.  A natural question raised, whether the augmented Lagrangian method can solve a general convex optimization problem with least constraint violation? This paper managed to answer this question within the general convex optimization framework. It is demonstrated that the dual of the convex  optimization problem with least constraint violation has an unbounded solution set, the optimality can be characterized by the augmented Lagrangian, and the augmented Lagrangian method is able to find an approximate solution, when the  least violated shift is in the domain of the subdifferential of the optimal value function.

There are many future works for the study of the optimization with the least constraint violation.
This paper only treats convex optimization problems. How about the augmented Lagrangian method for solving non-convex optimization problems with least constraint violation?  Even if the dual of the optimization problem with least constraint violation has a non-empty solution set, this problem is a typical convex optimization problem in which the conventional constraint qualifications (Slater or generalized Slater condition) are not satisfied. How to solve such kind of convex optimization problems?

\medskip \noindent
{\bf Acknowledgments.} The authors are very grateful to Prof. Ya-xiang Yuan for his long time guidance and encouragement and for Profs. Xinwei Liu and Zhongwen Chen for their useful discussions and comments.

\end{document}